\documentclass{article}

\PassOptionsToPackage{numbers}{natbib}
\usepackage[preprint]{neurips_2024}

\usepackage{microtype}
\usepackage{graphicx}
\usepackage{subfigure}
\usepackage{booktabs} 
\usepackage{nicefrac}
\usepackage{apxproof}  
\usepackage{hyperref}

\usepackage{amsmath}
\usepackage{amssymb}
\usepackage{mathtools}
\usepackage{amsthm}

\usepackage[capitalize,noabbrev]{cleveref}

\theoremstyle{plain}
\newtheoremrep{theorem}{Theorem}
\newtheoremrep{example}{Example}

\usepackage{todonotes}

\newcommand{\R}{\mathbb{R}}
\newcommand{\M}{\mathcal{M}}
\newcommand{\K}{\mathcal{K}}  
\newcommand{\C}{\mathcal{C}}  
\newcommand{\Q}{\mathcal{Q}}  
\newcommand{\Qr}{\mathcal{Q}_{r}}  
\newcommand{\Ksdp}{\mathcal{S}_{+}}

\newcommand{\Kexp}{\mathcal{E}}
\newcommand{\Kpow}{\mathcal{P}}

\newcommand{\e}{\mathbf{e}}  

\newcommand{\kgeq}{\succeq}
\newcommand{\kge}{\succ}

\DeclareMathOperator{\cl}{cl}  
\DeclareMathOperator{\intr}{int}  
\DeclareMathOperator{\relu}{ReLU}

\newcommand{\eigmin}{\lambda_{\min}}

\newcommand{\xbar}{\bar{x}}
\newcommand{\xhat}{\hat{x}}

\newcommand{\yhat}{\hat{y}}
\newcommand{\ybar}{\bar{y}}
\newcommand{\zhat}{\hat{z}}
\newcommand{\zlhat}{\hat{z}^{l}}
\newcommand{\zuhat}{\hat{z}^{u}}
\newcommand{\zl}{z^{l}}
\newcommand{\zu}{z^{u}}

\title{Dual Lagrangian Learning for Conic Optimization}

\author{%
  Mathieu Tanneau\\
  H. Milton Steward School of Industrial and Systems Engineering\\
  NSF AI Institute for Advances in Optimization\\
  Georgia Institute of Technology, Atlanta, GA, USA \\
  \texttt{mathieu.tanneau@isye.gatech.edu} \\
  \And
  Pascal Van Hentenryc\\
  H. Milton Steward School of Industrial and Systems Engineering\\
  NSF AI Institute for Advances in Optimization\\
  Georgia Institute of Technology, Atlanta, GA, USA \\
  \texttt{pascal.vanhentenryck@isye.gatech.edu} \\
}

\begin{document}

\maketitle

\begin{abstract}
This paper presents Dual Lagrangian Learning (DLL), a principled learning methodology for dual conic optimization proxies.
DLL leverages conic duality and the representation power of ML models to provide high-duality, dual-feasible solutions, and therefore valid Lagrangian dual bounds, for linear and nonlinear conic optimization problems.
The paper introduces a systematic dual completion procedure, differentiable conic projection layers, and a self-supervised learning framework based on Lagrangian duality.
It also provides closed-form dual completion formulae for broad classes of conic problems, which eliminate the need for costly implicit layers.
The effectiveness of DLL is demonstrated on linear and nonlinear conic optimization problems.
The proposed methodology significantly outperforms a state-of-the-art learning-based method, and achieves 1000x speedups over commercial interior-point solvers with optimality gaps under 0.5\% on average.
\end{abstract}

\section{Introduction}
\label{sec:introduction}

From power systems and manufacturing to supply chain management,
logistics and healthcare, optimization technology underlies most
aspects of the economy and society.  Over recent years, the
substantial achievements of Machine Learning (ML) have spurred
significant interest in combining the two methodologies.  This
integration has led to the development of new optimization algorithms
(and the revival of old ones) taylored to ML problems, as well as new
ML techniques for improving the resolution of hard optimization
problems \cite{Bengio2021_ML4CO}.  This paper focuses on the latter
(ML for optimization), specifically, the development of so-called
\emph{optimization proxies}, i.e., ML models that provide approximate
solutions to parametric optimization problems, see e.g.,
\cite{Kotary2021_E2EOPTLearningSurvey}.

In that context, considerable progress has been made in learning
\emph{primal} solutions for a broad range of problems, from linear to
discrete and nonlinear, non-convex optimization problems.
State-of-the-art methods can now predict high-quality, feasible or
close-to-feasible solutions for various applications
\cite{Kotary2021_E2EOPTLearningSurvey}.  This paper complements these
methods by learning \emph{dual} solutions which, in turn, certify the
(sub)optimality of learned primal solutions.  Despite the fundamental
role of duality in optimization, there is no dedicated framework for
dual optimization proxies, which have seldom received any attention in
the literature.  The paper addresses this gap by proposing, for the
first time, a principled learning methodology that combines conic
duality theory with Machine Learning. As a result, it becomes
possible, for a large class of optimization problems, to design a
primal proxy to deliver a high-quality primal solution and an
associated dual proxy to obtain a quality certificate.

\subsection{Contributions and outline}

The core contribution of the paper is the Dual Lagrangian Learning
(DLL) methodology for learning dual-feasible solutions for parametric
conic optimization problems.  DLL leverages conic duality to design a
self-supervised Lagrangian loss for training dual conic optimization
proxies.  In addition, the paper proposes a general dual conic
completion using differential conic projections and implicit
layers to guarantee dual feasibility, which yields stronger
guarantees than existing methods for constrained optimization
learning.  Furthermore, it presents closed-form analytical solutions
for conic projections, and for dual conic completion across broad classes of problems.
This eliminates the need for implicit layers in practice.
Finally, numerical results on linear and nonlinear conic problems demonstrate the effectiveness
of DLL, which a outperforms state-of-the-art learning baseline, and
yields significant speedups over interior-point solvers.

The rest of the paper is organized as follows.
Section \ref{sec:literature} presents the relevant literature.
Section \ref{sec:background} introduces notations and background material.
Section \ref{sec:methodology} presents the DLL methodology, which comprises the Lagrangian loss, dual completion strategy, and conic projections.
Section \ref{sec:exp} reports numerical results. Section \ref{sec:limitations} discusses the limitations of DLL, and Section \ref{sec:conclusion} concludes the paper.

\section{Related works}
\label{sec:literature}

    \paragraph{Constrained Optimization Learning}
        The vast majority of existing works on optimization proxies focuses on learning \emph{primal} solutions and, especially, on ensuring their feasibility.
        This includes, for instance, physics-informed training loss \cite{Fioretto2020_ACOPF_LG,Pan2020_deepopf,Donti2021_DC3}, mimicking the steps of an optimization algorithm \cite{Park2023_PDL,Qian2023_GNN4IPM,kotary2024_DeepAugmentedLagrangian}, using masking operations \cite{Bello2016_ML4TSP,Khalil2017_LearningCO}, or designing custom projections and feasibility layers \cite{Chen2023_E2EFeasibleProxies,Tordesillas2023_rayen}.
        The reader is referred to \cite{Kotary2021_E2EOPTLearningSurvey} for an extensive survey of constrained optimization learning.

        Only a handful of methods offer feasibility guarantees, and only for convex constraints; this latter point is to be expected since satisfying non-convex constraints is NP-hard in general.
        Implicit layers \cite{Agrawal2019_DifferentiableConvexLayers} have a high computational cost, and are therefore impractical unless closed-form solutions are available.
        DC3 \cite{Donti2021_DC3} uses equality completion and inequality correction, and is only guaranteed to converge for convex constraints and given enough correction steps.
        LOOP-LC \cite{Li2022_LOOP-LC} uses a gauge mapping to ensure feasibility for bounded polyhedral domains.
        RAYEN \cite{Tordesillas2023_rayen} and the similar work in \cite{konstantinov2024_NNHardConstraints} use a line search-based projection mechanism to handle convex constraints.
        All above methods employ equality completion, and the latter three \cite{Li2022_LOOP-LC,Tordesillas2023_rayen,konstantinov2024_NNHardConstraints} assume knowledge of a strictly feasible point, which is not always available.

    
    \paragraph{Dual Optimization Learning}

        To the authors' knowledge, dual predictions have received very little attention, with most works using them to warm-start an optimization algorithm.
        In \cite{Mak2023_ML-ADMM-OPF}, a primal-dual prediction is used to warm-start an ADMM algorithm, while  \cite{Kraul2023_MLDualCuttingStock} and \cite{Sugishita2023_MLWarmstartColGen} use a dual prediction as warm-start in a column-generation algorithm.
        These works consider specific applications, and do not provide dual feasibility guarantees.
        More recently, \cite{Park2023_PDL,kotary2024_DeepAugmentedLagrangian} attempt to mimic the (dual) steps of an augmented Lagrangian method, however with the goal of obtaining high-quality primal solutions.

        The first work to explicitly consider dual proxies, and to offer dual feasibility guarantees, is  \cite{Qiu2023_DualConicProxies_ACOPF}, which learns a dual proxy for a second-order cone relaxation of the AC Optimal Power Flow.
        Klamkin et al. \cite{klamkin2024dual} later introduce a dual interior-point learning algorithm to speed-up the training of dual proxies for bounded linear programming problems.
        In contrast, this paper proposes a general methodology for conic optimization problems, thus generalizing the approach in \cite{Qiu2023_DualConicProxies_ACOPF}, and provides more extensive theoretical results.
        The dual completion procedure used in \cite[Lemma 1]{klamkin2024dual} is a special case of the one proposed in this paper.

\section{Background}
\label{sec:background}

    This section introduces relevant notations and standard results on conic optimization and duality, which lay the basis for the proposed learning methodology.
    The reader is referred to \cite{ben2001lectures} for a thorough overview of conic optimization.

\subsection{Notations}
\label{sec:background:notations}

    Unless specified otherwise, the Euclidean norm of a vector $x \in \R^{n}$ is denoted by $\|x\| = \sqrt{x^{\top}x}$.
    The positive and negative part of $x \in \R$ are denoted by $x^{+} = \max(0, x)$ and $x^{-} = \max(0, -x)$.
    The identity matrix of order $n$ is denoted by $I_{n}$, and $\e$ denotes the vector of all ones.
    The smallest eigenvalue of a real symmetric matrix $X$ is $\eigmin(X)$.

    Given a set $\mathcal{X} \subseteq \R^{n}$, the \emph{interior} and \emph{closure} of $\mathcal{X}$ are denoted by $\intr \mathcal{X}$ and by $\cl \mathcal{X},$ respectively.
    The Euclidean projection onto convex set $\mathcal{C}$ is denoted by $\Pi_{\mathcal{C}}$, where
    \begin{align}
        \label{eq:projection}
         \displaystyle \Pi_{\mathcal{C}}(\bar{x}) = \arg \min_{x \in \mathcal{C}} \|x - \bar{x}\|^{2} \ .
    \end{align}
    The set $\K \,{\subseteq} \, \R^{n}$ is a \emph{cone} if $x \in \K, \lambda \geq 0 \Rightarrow \lambda x \in \K$.
    The \emph{dual cone} of $\K$ is
    \begin{align}
        \K^{*} = \{y \in \R^{n}: y^{\top}x \geq 0, \forall x \in \K\},
    \end{align}
    whose negative $\K^{\circ} = -\K^{*}$ is the \emph{polar} cone of $\K$.
    A cone $\K$ is self-dual if $\K \, {=} \, \K^{*}$, and it is \emph{pointed} if $\K \cap (-\K) \, {=} \, \{0\}$.
    All cones considered in the paper are \emph{proper} cones, i.e., closed, convex, pointed cones with non-empty interior.
    A proper cone $\K$ defines \emph{conic inequalities} $\kgeq_{\K}$ and $\kge_{\K}$ as 
    \begin{subequations}
        \label{eq:def:conic_inequalities}
        \begin{align}
            \forall (x, y) \in \R^{n} {\times} \R^{n},& \  x \kgeq_{\K} y \, \Leftrightarrow \, x - y \in \K,\\
            \forall (x, y) \in \R^{n} {\times} \R^{n},& \  x \kge_{\K}  y \, \Leftrightarrow \, x - y \in \intr \K.
        \end{align}
    \end{subequations}
    
\subsection{Conic optimization}
\label{sec:background:conic}

    Consider a (convex) conic optimization problem of the form
    \begin{align}
        \label{eq:background:conic:primal}
        \min_{x} \quad 
        \left\{
            c^{\top} x
        \ \middle| \
            A x \succeq_{\K} b
        \right\},
    \end{align}
    where $A \, {\in} \, \R^{m \times n}, b \, {\in} \, \R^{m}, c \, {\in} \, \R^{n}$, and $\K$ is a proper cone.
    All convex optimization problems can be formulated in conic form.
    A desirable property of conic formulations is that is enables the use of principled conic duality theory \cite{ben2001lectures}.
    Namely, the conic dual problem reads
    \begin{align}
        \label{eq:background:conic:dual}
        \max_{y} \quad 
        \left\{
            b^{\top} y
        \ \middle| \
            A^{\top} y = c, y \in \K^{*}
        \right\}.
    \end{align}
    The dual problem \eqref{eq:background:conic:dual} is a conic problem, and the dual of \eqref{eq:background:conic:dual} is  \eqref{eq:background:conic:primal}.
    Weak conic duality always holds, i.e., any dual-feasible solution provides a valid lower bound on the optimal value of \eqref{eq:background:conic:primal}, and vice-versa.
    When strong conic duality holds, e.g., under Slater's condition, both primal/dual problems have the same optimal value and a primal-dual optimal solution exists \cite{ben2001lectures}.
    
    Conic optimization encompasses broad classes of problems such linear and semi-definite programming.
    Most real-life convex optimization problems can be represented in conic form using only a small number of cones \cite{Lubin2016_ExtendedFormulationsMICONV}, which are supported by off-the-shelf solvers such as Mosek, ECOS, or SCS.
    These so-called ``standard" cones comprise the non-negative orthant $\R_{+}$, the second-order cone $\Q$ and rotated second-order cone $\Qr$, the positive semi-definite cone $\Ksdp$, the power cone $\Kpow$ and the exponential cone $\Kexp$; see Appendix \ref{app:cones} for algebraic definitions.

\section{Dual Lagrangian Learning (DLL)}
\label{sec:methodology}

\begin{figure}
    \centering
    \includegraphics[width=\columnwidth]{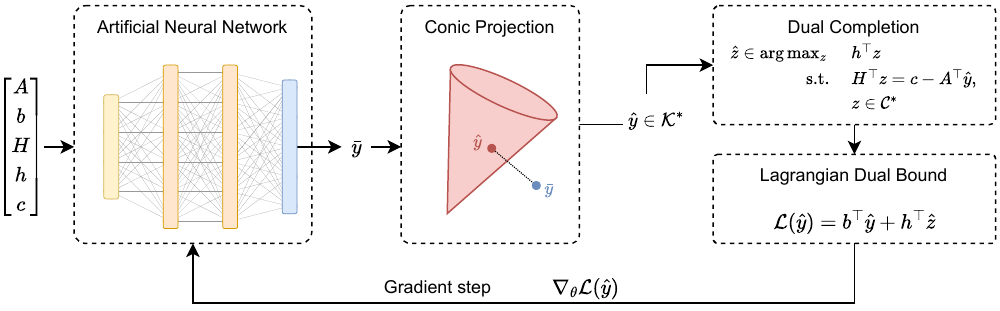}
    \caption{Illustration of the proposed DLL scheme.
    Given input data $(A, b, H, h, c)$, a neural network first predicts $\ybar \in \R^{n}$.
    Next, a conic projection layer computes a conic-feasible $\yhat \in \K^{*}$, which is then completed into a full dual-feasible solution $(\yhat, \zhat)$.
    The model is trained in a self-supervised fashion, by updating the weights $\theta$ to maximize the Lagrangian dual bound $\mathcal{L}(\yhat)$.}
    \label{fig:DLL}
\end{figure}

This section presents the {\em Dual Lagrangian Learning} (DLL) methodology, illustrated in Figure \ref{fig:DLL}, for learning dual solutions for conic optimization problems.
DLL combines the representation power of artificial neural networks (or, more generally, any differentiable program), with conic duality theory, thus providing valid Lagrangian dual bounds for general conic optimization problems.
{\em To the best of the authors' knowledge, this paper is the first to propose a principled self-supervised framework with dual guarantees for general conic optimization problems.}

DLL exploits three fundamental building blocks: (1) a dual conic
completion procedure that provides dual-feasible solutions and, hence,
valid Lagrangian dual bounds; (2) fast and differentiable conic
projection layers; and (3) a self-supervised learning algorithm that
emulates the steps of a dual Lagrangian ascent algorithm.

\subsection{Dual Conic Completion}
\label{sec:methodology:completion}

    Consider a conic optimization problem in primal-dual form\\
    \noindent
    \begin{minipage}{0.45\textwidth}
    \begin{subequations}
        \label{eq:completion:conic:primal}
        \begin{align}
            \min_{x} \quad & c^{\top} x\\
            \text{s.t.} \quad 
            & \label{eq:completion:conic:primal:hard_constraints}
                A x \kgeq_{\K} b,\\
            & \label{eq:completion:conic:primal:bounding_constraints}
                H x \kgeq_{\C} h,
        \end{align}
    \end{subequations}
    \end{minipage}
    \hfill
    \begin{minipage}{0.45\textwidth}
    \begin{subequations}
        \label{eq:completion:conic:dual}
        \begin{align}
        \label{eq:completion:conic:dual:obj}
            \max_{y, z} \quad & b^{\top} y + h^{\top}z\\
            s.t. \quad 
            & \label{eq:completion:conic:dual:con:eq}
            A^{\top} y + H^{\top}z = c,\\
            & \label{eq:completion:conic:dual:con:conic}
                y \in \K^{*}, z \in \C^{*}.
        \end{align}
    \end{subequations}
    \end{minipage}

    \noindent
    where $y \, {\in} \, \K^{*}$ and $z \, {\in} \, \C^{*}$ are the dual variables associated to constraints \eqref{eq:completion:conic:primal:hard_constraints} and \eqref{eq:completion:conic:primal:bounding_constraints}, respectively.
    The proposed dual conic completion, outlined in Theorems \ref{thm:completion:existence} and \ref{thm:completion:optimal} below, takes as input $\yhat \, {\in} \, \K^{*}$, and recovers $\zhat \, {\in} \, \C^{*}$ such that $(\yhat, \zhat)$ is feasible for \eqref{eq:completion:conic:dual}.   
    The initial assumption that $\yhat \, {\in} \, \K^{*}$ can be enforced through a projection step, which will be described in Section \ref{sec:methodology:projection}.
    
    \begin{theorem}[Dual conic completion]
        \label{thm:completion:existence}
        Assume that $\forall \yhat \in \K^{*}, \exists x : Hx \kge_{\C} h$ and the problem
        \begin{align}
        \label{eq:completion:lagrange:primal}
            \min_{x} \quad 
            \left\{ 
                c^{\top}x + (b - Ax)^{\top} \yhat
            \ \middle| \
                Hx \kgeq_{\C} h
            \right\}
        \end{align}
        is bounded.
        Then,
        $
            \forall \yhat \in \K^{*},
            \exists \zhat \, {\in} \, \C^{*}:
            A^{\top} \yhat + H^{\top} \zhat = c
        $,
        i.e., $(\yhat, \zhat)$ is feasible for \eqref{eq:completion:conic:dual}.
    \end{theorem}
    \begin{proof}
        Let $\yhat \in \K^{*}$, and recall that \eqref{eq:completion:lagrange:primal} is bounded and strictly feasible.
        By strong duality, its dual
        \begin{align}
        \label{eq:completion:lagrange:dual}
            \max_{z} \quad 
            \left\{
                h^{\top}z + b^{\top} \yhat
            \ \middle| \ 
                H^{\top}z = c - A^{\top} \yhat,
                \,
                z \in \C^{*}
            \right\}
        \end{align}
        is solvable \cite{ben2001lectures}.
        Therefore, there exists a feasible solution $\zhat$ for \eqref{eq:completion:lagrange:dual}.
        By construction, $\zhat \in \C^{*}$ and $A^{\top} \yhat + H^{\top} \zhat = c$, hence $(\yhat, \zhat)$ is feasible for \eqref{eq:completion:conic:dual}.
    \end{proof}

    \begin{theorem}[Optimal dual completion]
        \label{thm:completion:optimal}
        Let $\yhat \in \K^{*}$, and let $\zhat$ be dual-optimal for \eqref{eq:completion:lagrange:primal}.\\
        Then, $\mathcal{L}(\hat{y}, \zhat) = b^{\top} \yhat + h^{\top} \zhat$ is a valid dual bound on the optimal value of \eqref{eq:completion:conic:primal}, and $\mathcal{L}(\yhat, \zhat)$ is the strongest dual bound that can be obtained after fixing $y = \yhat$ in \eqref{eq:completion:conic:dual}.
    \end{theorem}
    \begin{proof}
        First, recall that $\zhat$ exists by strong conic duality; see proof of Theorem \ref{thm:completion:existence}.
        Furthermore, $(\yhat, \zhat)$ is feasible for \eqref{eq:completion:conic:dual} by construction.
        Thus, by weak duality, the Lagrangian bound $\mathcal{L}(\yhat) = b^{\top} \yhat + h^{\top} \zhat$ is a valid dual bound on the optimal value of \eqref{eq:completion:conic:primal}.
        Finally, fixing $y = \yhat$ in \eqref{eq:completion:conic:dual} yields
        \begin{align}
            \label{eq:completion:dual:fixed_y}
            \max_{y, z} \quad
            \left\{
                \eqref{eq:completion:conic:dual:obj}
            \ \middle| \
                \eqref{eq:completion:conic:dual:con:eq},
                \eqref{eq:completion:conic:dual:con:conic},
                y = \yhat
            \right\},
        \end{align}
        which is equivalent to \eqref{eq:completion:lagrange:dual}.
        Hence, its optimal value is $b^{\top}\yhat {+} h^{\top}\zhat \, {=} \, \mathcal{L}(\yhat, \zhat)$ by definition of $\zhat$.
    \end{proof}

It is important to note the theoretical differences between the
proposed dual completion, and applying a generic method, e.g., DC3
\cite{Donti2021_DC3}, LOOP-LC \cite{Li2022_LOOP-LC} or RAYEN
\cite{Tordesillas2023_rayen}, to the dual problem
\eqref{eq:completion:conic:dual}.  First, LOOP-LC is not applicable
here, because it only handles linear constraints and requires a
compact feasible set, which is not the case in general for
\eqref{eq:completion:conic:dual}. Second, unlike RAYEN, Theorem
\ref{thm:completion:existence} does not require an explicit
characterization of the affine hull of the (dual) feasible set, nor
does it assume knowledge of a strictly feasible point.  In fact,
Theorem \ref{thm:completion:existence} applies even if the feasible
set of \eqref{eq:completion:conic:dual} has an empty interior.  Third,
the proposed dual completion enforces both linear equality constraints
\eqref{eq:completion:conic:dual:con:eq} and conic constraints
\eqref{eq:completion:conic:dual:con:conic}.  In contrast, the equality
completion schemes used in DC3 and RAYEN enforce equality constraints
but need an additional mechanism to handle inequality constraints.
Fourth, the optimal completion outlined in Theorem
\ref{thm:completion:optimal} provides guarantees on the strength of
the Lagrangian dual bound $\mathcal{L}(\yhat, \zhat)$.  This is a
major difference with DC3 and RAYEN, whose correction mechanism does
not provide any guarantee of solution quality.  {\em Overall, the
fundamental difference between generic methods and the proposed
optimal dual completion, is that the former only exploit dual
feasibility constraints
\eqref{eq:completion:conic:dual:con:eq}--\eqref{eq:completion:conic:dual:con:conic},
whereas DLL also exploits (dual) optimality conditions, thus
providing additional guarantees.}

    Another desirable property of the proposed dual completion procedure, is that it does not require the user to formulate the dual problem \eqref{eq:completion:conic:dual} explicitly, as would be the case for DC3 or RAYEN.
    Instead, the user only needs to identify a set of \emph{primal} constraints that satisfy the conditions of Theorem \ref{thm:completion:existence}.
    For instance, it suffices to identify constraints that bound the set of primal-feasible solutions.
    This is advantageous because practitioners typically work with primal problems rather than their dual.
    The optimal dual completion can then be implemented via an implicit optimization layer.
    Thereby, in a forward pass, $\zhat$ is computed by solving the primal-dual pair \eqref{eq:completion:lagrange:primal}--\eqref{eq:completion:lagrange:dual} and, in a backward pass, gradient information is obtained via the implicit function theorem \cite{Agrawal2019_DifferentiatingConeProgram}.

    The main limitations of implicit layers are their numerical instability and their computational cost, both in the forward and backward passes.
    To eliminate these issues, closed-form analytical solutions are presented next for broad classes of conic optimization problems; other examples are presented in the numerical experiments of Section \ref{sec:exp}.

    \begin{example}[Bounded variables]
        \label{ex:bounded_variables}
        Consider a conic optimization problem with bounded variables
        \begin{align}
            \min_{x} \quad
            \left\{ 
                c^{\top}x 
            \ \middle| \
                Ax \kgeq_{\K} b, l \leq x \leq u
            \right\}
        \end{align}
        where $l < u$ are finite lower and upper bounds on all variables $x$.
        The dual problem is
        \begin{align}
            \min_{y, \zl, \zu} \quad
            \left\{ 
                b^{\top}y + l^{\top} \zl - u^{\top} \zu
            \ \middle| \
                A^{\top} y + \zl - \zu = c, y \in \K^{*}, \zl \geq 0, \zu \geq 0
            \right\}
        \end{align}
        and the optimal dual completion is $\zlhat = |c - A^{\top} \yhat|^{+}, \zuhat = |c - A^{\top} \yhat|^{-}$.
    \end{example}
    \begin{proof}
        Let $\yhat \in \K^{*}$ be fixed.
        Fixing $y = \yhat$ in the dual problem yields
        \begin{align}
            \max_{\zl, \zu} \quad & l^{\top} \zl - u^{\top} \zu\\
            \text{s.t.} \quad
            & \zl - \zu = c - A^{\top} \yhat\\
            & \zl, \zu \geq 0.
        \end{align}
        Eliminating $\zl = \zu + (c - A^{\top}\yhat)$, the problem becomes
        \begin{align}
            \max_{\zu} \quad & (l - u)^{\top} \zu + l^{\top}(c - A^{\top}\yhat)\\
            \text{s.t.} \quad
            & \zu \geq -(c - A^{\top} \yhat)\\
            & \zu \geq 0.
        \end{align}
        Since $l < u$, i.e., the objective coefficient of $\zu$ is negative, and the problem is a maximization problem, it follows that $\zu$ must be as small as possible in any optimal solution.
        Hence, at the optimum, $\zuhat = \max(0, -(c - A^{\top}y)) = |c - A^{\top} \yhat|^{-}$, and $\zlhat = |c - A^{\top}\yhat|^{+}$.
    \end{proof}

    The assumption in Example \ref{ex:bounded_variables} that all variables have finite bounds holds in most --if not all-- real-life settings, where decision variables are physical quantities (e.g. budgets or production levels) that are naturally bounded.
    The resulting completion procedure is a generalization of that used in \cite{klamkin2024dual} for linear programming (LP) problems.
    
    \begin{example}[Trust region]
        \label{ex:trust_region}
        Consider the trust region problem \cite{nocedal1999numerical}
        \begin{align}
            \min_{x} \quad
            \left\{
                c^{\top}x 
            \ \middle| \
                Ax \kgeq_{\K} b, \| x \| \leq r
            \right\}
        \end{align}
        where $r {\geq} 0$, $\| {\cdot} \|$ is a norm, and $\| x \| {\leq} r \Leftrightarrow (r, x) \, {\in} \, \C \, {=} \, \{(t, x) \, | \, t {\geq} \| x \|\}$.
        The dual problem is
        \begin{align}
            \max_{y, z_{0}, z} \quad 
            \left\{
                b^{\top} y - r z_{0}
            \ \middle| \
                A^{\top} y + z = c, 
                y \in \K^{*}, 
                (z_{0}, z) \in \C^{*}
            \right\}
        \end{align}
        where $\| {\cdot} \|_{*}$ is the dual norm and $\C^{*} \, {=} \, \{ (t, x) \, | \, t {\geq} \|x \|_{*}\}$ \cite{boyd2004convex}.
        The optimal dual completion is $\zhat = c - A^{\top} \yhat$, $\zhat_{0} = \| \zhat \|_{*}$.
    \end{example}
    \begin{proof}
        The relation $\zhat = c - A^{\top} \yhat$ is immediate from the dual equality constraint $A^{\top}y + z = c$.
        Next, observe that $z_{0}$ appears only in the constraint $(z_{0}, z) \in \C^{*}$, and has negative objective coeffient.
        Hence, $z_{0}$ must be as small as possible in any optimal solution.
        This yields $\zhat_{0} = \| \zhat \|_{*}$.
    \end{proof}
    
    \begin{example}[Convex quadratic objective]
        Consider the convex quadratic conic problem
        \begin{align}
            \min_{x} \quad
            \left\{ 
                \nicefrac{1}{2} \times
                x^{\top}Qx + c^{\top}x
            \ \middle| \ 
                Ax \kgeq_{\K} b
            \right\},
        \end{align}
        where $Q = F^{\top} F$ is positive definite.
        The problem can be formulated as the conic problem
        \begin{align}
            \min_{x} \quad
            \left\{ 
                q + c^{\top}x
            \ \middle| \ 
                Ax \kgeq_{\K} b,
                (1, q, Fx) \in \Qr^{2+n}
            \right\}
        \end{align}
        whose dual is
        \begin{align}
            \max_{y, z_{0}, z} \quad
            \left\{ 
                b^{\top}y - z_{0}
            \ \middle| \ 
                A^{\top}y + F^{\top}z = c,
                (1, z_{0}, z) \in \Qr^{2+n}
            \right\}.
        \end{align}
        The optimal dual completion is $\zhat = F^{-\top}(c - A^{\top} \yhat), \zhat_{0} = \nicefrac{1}{2} \| \zhat \|^{2}_{2}$.
    \end{example}
    \begin{proof}
        The proof uses the same argument as the proof of Example \ref{ex:trust_region}.
        Namely, $\zhat = F^{-\top}(c - A^{\top} \yhat)$ is immediate from the dual equality constraint $A^{\top}y + F^{\top}z = c$.
        Note that $F^{\top}$ is non-singular because $Q$ is positive definite.
        Finally, $z_{0}$ has negative objective coefficient, and only appears in the conic constraint $(1, z_{0}, z) \in \Qr^{2+n}$.
        Therefore, at the optimum, one must have $2 \zhat_{0} = \| \zhat \|_{2}^{2}$, which concludes the proof.
    \end{proof}

\subsection{Conic Projections}
\label{sec:methodology:projection}

    The second building block of DLL are differentiable conic projection layers.
    Note that DLL only requires a valid projection onto $\K^{*}$, which need not be the Euclidean projection $\Pi_{\K^{*}}$.
    Indeed, the latter may be computationally expensive and cumbersome to differentiate.
    For completeness, the paper presents Euclidean and non-Euclidean projection operators, where the latter are simple to implement, computationally fast, and differentiable almost everywhere.
    Closed-form formulae are presented for each standard cone in Appendix \ref{app:cones}, and an overview is presented in Table \ref{tab:conic_projections}.

    \begin{table}[!t]
        \centering
        \caption{Overview of conic projections for standard cones}
        \label{tab:conic_projections}
        \begin{tabular}{cccc}
        \toprule
            Cone & Definition & Euclidean projection & Radial projection\\
        \midrule
            $\R_{+}$ & Appendix \ref{app:cones:R+} & \eqref{eq:cones:R+:proj:euclidean} & \eqref{eq:cones:R+:proj:euclidean}\\
            $\Q$     & Appendix \ref{app:cones:SOC} & \eqref{eq:cones:SOC:proj:euclidean} & \eqref{eq:cones:SOC:proj:radial} \\
            $\Ksdp$  & Appendix \ref{app:cones:SDP} & \eqref{eq:cones:SDP:proj:euclidean} & \eqref{eq:cones:SDP:proj:radial} \\
            $\Kexp$  & Appendix \ref{app:cones:exp} & no closed form & \eqref{eq:cones:exp:proj:radial:primal} and \eqref{eq:cones:exp:proj:radial:dual} \\
            $\Kpow$  & Appendix \ref{app:cones:pow} & no closed form & \eqref{eq:cones:pow:proj:radial} \\
        \bottomrule
        \end{tabular}
    \end{table}

    \subsubsection{Euclidean projection}

        Let $\K$ be a proper cone, and $\xbar \in \R^{n}$.
        By Moreau's decomposition \cite{Parikh2014_ProximalAlgorithms},
        \begin{align}
            \label{eq:projection:moreau}
            \xbar &= \Pi_{\K}(\xbar) + \Pi_{\K^{\circ}}(\xbar),
        \end{align}
        which is a reformulation of the KKT conditions of the projection problem \eqref{eq:projection}, i.e.,
        \begin{align}
            \xbar = p - q, \quad p \in \K, \quad q \in \K^{*}, \quad p^{\top}q = 0.
        \end{align}
        It then follows that $\Pi_{\K^{*}}(\xbar) = - \Pi_{\K^{\circ}}(-\xbar)$, by invariance of Moreau's decomposition under orthogonal transformations.
        Thus, it is sufficient to know how to project onto $\K$ to be able to project onto $\K^{*}$ and $\K^{\circ}$.
        Furthermore, \eqref{eq:projection:moreau} shows that $\Pi_{\K}$ is identically zero on the polar cone $\K^{\circ}$.
        In a machine learning context, this may cause gradient vanishing issues and slow down training.

    \subsubsection{Radial projection}

        Given an interior ray $\rho \kge_{\K} 0$, the radial projection operator $\Pi_{\K}^{\rho}$ is defined as
        \begin{align}
             \Pi_{\K}^{\rho}(\xbar) = \xbar + \lambda \rho \quad \text{where} \quad \lambda = \min_{\lambda \geq 0} \{\lambda \, | \, \xbar + \lambda \rho \in \K\}.
        \end{align}
        The name stems from the fact that $\Pi_{\K}^{\rho}$ traces ray $\rho$ from $\xbar$ until $\K$ is reached.
        Unlike the Euclidean projection, it requires an interior ray, which however only needs to be determined once \emph{per cone}.
        The radial projection can then be computed, in general, via a line search on $\lambda$ or via an implicit layer.
        Closed-form formulae for standard cones and their duals are presented in Appendix \ref{app:cones}.

\subsection{Self-Supervised Dual Lagrangian Training}
\label{sec:methodology:training}

    The third building block of DLL is a self-supervised learning framework for training dual conic optimization proxies.
    In all that follows, let $\xi = (A, b, H, h, c)$ denote the data of an instance \eqref{eq:completion:conic:primal}, and assume a distribution of instances $\xi \sim \Xi$.
    Next, let $\M_{\theta}$ be a differentiable program parametrized by $\theta$, e.g., an artificial neural network, which takes as input $\xi$ and outputs a dual-feasible solution $(\yhat, \zhat)$.
    Recall that dual feasibility of $(\yhat, \zhat)$ can be enforced by combining the dual conic projection presented in Section \ref{sec:methodology:projection}, and the optimal dual completion outlined in Theorem \ref{thm:completion:optimal}.
    
    The proposed self-supervised dual lagrangian training is formulated as
    \begin{subequations}
        \label{eq:training_self_supervised}
        \begin{align}
            \max_{\theta} \quad 
                    & \mathbb{E}_{\xi \sim \Xi} \left[
                        \mathcal{L}(\yhat, \zhat, \xi) \right]\\
            \text{s.t.} \quad
                    & (\yhat, \zhat) = \M_{\theta}(\xi),
        \end{align}
    \end{subequations}
    where $\mathcal{L}(\yhat, \zhat, \xi) \, {=} \, b^{\top} \yhat \, {+} \, h^{\top} \zhat$ is the Lagrangian dual bound obtained from $(\yhat, \zhat)$ by weak duality.
    Thereby, the training problem \eqref{eq:training_self_supervised} seeks the value of $\theta$ that maximizes the expected Lagrangian dual bound over the distribution of instances $\Xi$, effectively mimicking the steps of a (sub)gradient algorithm.
    Note that, instead of updating $(\yhat, \zhat)$ directly, the training procedure computes a (sub)gradient $\partial_{\theta} \mathcal{L}(\yhat, \zhat, \xi)$ to update $\theta$, and then obtains a new prediction $(\yhat, \zhat)$ through $\M_{\theta}$.
    Also note that formulation \eqref{eq:training_self_supervised} does not required labeled data, i.e., it does not require pre-computed dual-optimal solutions.
    Furthermore, it applies to any architecture that guarantees dual feasibility of $(\yhat, \zhat)$, i.e., it does not assume any specific projection nor completion procedure.

\section{Numerical experiments}
\label{sec:exp}

This section presents numerical experiments on linear and nonlinear optimization problems; detailed problem formulations, model architectures, and other experiment settings, are reported in Appendix \ref{sec:app:experiment_details}.
The proposed DLL methodology is evaluated against applying DC3 to the dual problem \eqref{eq:completion:conic:dual} as a baseline.
Thereby, linear equality constraints \eqref{eq:completion:conic:dual:con:eq} and conic inequality constraints \eqref{eq:completion:conic:dual:con:conic} are handled by DC3's equality completion and inequality correction mechanisms, respectively.
The two approaches (DLL and DC3) are evaluated in terms of dual optimality gap and training/inference time.
The dual optimality gap is defined as $(\mathcal{L}^{*} - \mathcal{L}(\yhat, \zhat)) / \mathcal{L}^{*}$, where $\mathcal{L}^{*}$ is the optimal value obtained from a state-of-the-art interior-point solver.

\subsection{Linear Programming Problems}
\label{sec:mdk}

    \subsubsection{Problem formulation and dual completion}
    \label{sec:mdk:formulation}

    The first set of experiments considers continuous relaxations of multi-dimensional knapsack problems \cite{Freville1994_MDKnapsackPreprocessing,Freville2004_MultiDimensionalKnapsackOverview}, which are of the form
    \begin{align}
    \label{eq:multi-knapsack}
        \min_{x} \quad
        \left\{
            -p^{\top}x
        \ \middle| \ 
            Wx \leq b,
            x \in [0, 1]^{n}
        \right\}
    \end{align}
    where $p \in \R_{+}^{n}$, $W \in \R_{+}^{m \times n}$, and $b \in \R_{+}^{m}$.
    The dual problem reads
    \begin{align}
        \label{eq:multi-knapsack:dual}
        \max_{y, \zl, \zu} \quad
        \left\{
            b^{\top} y - \e^{\top} \zu
        \ \middle| \
            W^{\top}y + \zl - \zu = -p,
            y \leq 0, \zl \geq 0, \zu \geq 0,
        \right\}
    \end{align}
    where $y \in \R^{m}$ and $\zl, \zu \in \R^{n}$.
    Since variables $x$ is bounded, the closed-form completion presented in Example \ref{ex:bounded_variables} applies.
    Namely, $\zlhat = | {-}p {-} W^{\top} \yhat|^{+}$ and $\zuhat = |{-}p {-} W^{\top} \yhat|^{-}$, where $\yhat \in \R_{-}^{m}$.

    \subsubsection{Numerical results}
    \label{sec:mdk:results}

    \begin{table}[!t]
        \centering
        \caption{Comparison of optimality gaps on linear programming instances.}
        \label{tab:res:MDK:gaps}
        \small
        \begin{tabular}{ccrrrrrrrr}
            \toprule
                &&
                & \multicolumn{3}{c}{DC3}
                & \multicolumn{3}{c}{DLL}\\
            \cmidrule(lr){4-6}
            \cmidrule(lr){7-9}
            \multicolumn{1}{c}{$m$} & \multicolumn{1}{c}{$n$}
                & Opt val$^*$
                & avg. & std & max 
                & avg. & std & max\\
            \midrule
               5 &  100 &  14811.9 &  19.58 &   1.86 &  41.42 &   \textbf{0.36} &   0.20 &   1.36 \\
                 &  200 &  29660.4 &  20.58 &   1.41 &  49.47 &   \textbf{0.18} &   0.10 &   0.84 \\
                 &  500 &  74267.0 &  33.70 &   1.29 &  41.54 &   \textbf{0.07} &   0.04 &   0.30 \\
            \midrule
              10 &  100 &  14675.8 &  41.85 &   2.51 &  69.58 &   \textbf{0.68} &   0.25 &   2.15 \\
                 &  200 &  29450.7 &  36.88 &   2.28 & 100.90 &   \textbf{0.34} &   0.13 &   0.96 \\
                 &  500 &  73777.5 & 100.04 &   3.38 & 104.00 &   \textbf{0.14} &   0.06 &   0.46 \\
            \midrule
              30 &  100 &  14441.5 & 159.49 &   5.54 & 166.31 &   \textbf{1.93} &   0.37 &   3.31 \\
                 &  200 &  29156.1 & 255.24 &   8.42 & 259.25 &   \textbf{0.96} &   0.20 &   1.83 \\
                 &  500 &  73314.3 & 274.78 &   7.91 & 277.40 &   \textbf{0.38} &   0.09 &   0.75 \\
             \bottomrule
        \end{tabular}\\
        \footnotesize{All gaps are in \%; best values are in bold. $^*$Mean optimal value on test set; obtained with Gurobi.}
    \end{table}
    
    Table \ref{tab:res:MDK:gaps} reports, for each combination of $m, n$: the average optimal value obtained by Gurobi (Opt val), as well as the average (avg), standard-deviation (std) and maximum (max) optimality gaps achieved by DC3 and DLL on the test set.
    First, DLL significantly outperforms DC3, with average gaps ranging from 0.07\% to 1.93\%, compared with 19.58\%--274.78\% for DC3, an improvement of about two orders of magnitude.
    A similar behavior is observed for maximum optimality gaps.
    The rest of the analysis thus focuses on DLL.
    Second, an interesting trend can be identified: optimality gaps tend to increase with $m$ and decrease with $n$.
    This effect may be explained by the fact that increasing $m$ increases the output dimension of the FCNN; larger output dimensions are typically harder to predict.
    In addition, a larger $n$ likely provides a smoothing effect on the dual, whose solution becomes easier to predict.
    The reader is referred to \cite{Freville2004_MultiDimensionalKnapsackOverview} for probabilistic results on properties of multi-knapsack problems.

        \begin{table}[!t]
            \centering
            \caption{Computing time statistics for linear programming instances}
            \label{tab:res:MDK:timing}
            \small
            \begin{tabular}{rrrrrrrrr}
                \toprule
                 $m$ & $n$ 
                    & \multicolumn{1}{c}{Gurobi$^{\dagger}$} 
                    & \multicolumn{1}{c}{DC3$^{\ddagger}$} 
                    & \multicolumn{1}{c}{DLL$^{\ddagger}$}\\
\midrule
   5 &  100 &      2.8 CPU.s &    2.1 GPU.ms &    0.3 GPU.ms\\
     &  200 &      4.1 CPU.s &    4.0 GPU.ms &    0.7 GPU.ms\\
     &  500 &      6.6 CPU.s &   13.2 GPU.ms &    3.0 GPU.ms\\
\midrule
  10 &  100 &      3.7 CPU.s &    2.3 GPU.ms &    0.4 GPU.ms\\
     &  200 &      6.1 CPU.s &    4.9 GPU.ms &    1.1 GPU.ms\\
     &  500 &     11.9 CPU.s &   17.2 GPU.ms &    4.7 GPU.ms\\
\midrule
  30 &  100 &     14.0 CPU.s &    4.6 GPU.ms &    0.9 GPU.ms\\
     &  200 &     21.3 CPU.s &   10.4 GPU.ms &    2.5 GPU.ms\\
     &  500 &     40.0 CPU.s &   39.5 GPU.ms &   13.6 GPU.ms\\
                \bottomrule
            \end{tabular}\\
            {\footnotesize $^{\dagger}$Time to solve all instances in the test set, using one CPU core. $^{\ddagger}$Time to run inference on all instances in the test set, using one V100 GPU.}
        \end{table}

    Next, Table \ref{tab:res:MDK:timing} reports computing time statistics for Gurobi, DC3 and DLL.
    Namely, the table reports, for each combination of $m, n$, the time it takes to execute each method on all instances in the test set.
    First, DLL is 3-10x faster than DC3, which is caused by DC3's larger output dimension ($m{+}n$, compared to $m$ for DLL), and its correction steps.
    Furthermore, unsurprisingly, both DC3 and DLL yield substantial speedups compared to Gurobi, of about 3 orders of magnitude.
    Note however that Gurobi's timings could be improved given additional CPU cores, although both ML-based methods remain significantly faster using a single GPU.

\subsection{Nonlinear Production and Inventory Planning Problems}
\label{sec:exp:rcprod}
        
    \subsubsection{Problem formulation and dual completion}
    \label{sec:rcprod:formulation}

        The second set of experiments considers the nonlinear resource-constrained production and inventory planning problem \cite{Ziegler1982_ProductionPlanning,MosekModelingCookbook}.
        In primal-dual form, the problem reads\\
        \begin{minipage}[t]{0.45\textwidth}
        \small
        \begin{subequations}
        \label{eq:rcprod:conic:primal}
        \begin{align}
            \min_{x, t} \quad
            & d^{\top}x + f^{\top}t\\
            \text{s.t.} \quad
            & \label{eq:rcprod:conic:primal:resource}
                r^{\top}x \leq b,\\
            & \label{eq:rcprod:conic:primal:conic}
                (x_{j}, t_{j}, \sqrt{2}) \in \Qr^{3}, \ j {=} 1,...,n
        \end{align}
        \end{subequations}
        \end{minipage}
        \hfill
        \begin{minipage}[t]{0.53\textwidth}
        \small
        \begin{subequations}
            \label{eq:rcprod:conic:dual}
            \begin{align}
                \max_{y, \pi, \tau, \sigma} \quad
                & b y - \sqrt{2} \e^{\top} \sigma_{j}\\
                s.t. \quad 
                & \label{eq:rcprod:conic:dual:pi}
                    r y + \pi = d,\\
                & \label{eq:rcprod:conic:dual:tau}
                    \tau = f,\\
                & \label{eq:rcprod:conic:dual:y}
                    y \leq 0,\\
                & \label{eq:rcprod:conic:dual:rsoc}
                    (\pi_{j}, \tau_{j}, \sigma_{j}) \in \Qr^{3}, \ j{=}1, ..., n
            \end{align}
        \end{subequations}
        \end{minipage}
        
        \noindent
        where $r, d, f \in \R^{n}$ are positive vectors, and $b > 0$.
        Primal variables are $x, t \in \R^{n}$, and the dual variables associated to constraints \eqref{eq:rcprod:conic:primal:resource} and \eqref{eq:rcprod:conic:primal:conic} are $y \in \R_{-}$, and $\pi, \sigma, \tau \in \R^{n}$, respectively.
        
        Note that \eqref{eq:rcprod:conic:primal:conic} implies $x, t \geq 0$. 
        Next, let $y \leq 0$ be fixed, and consider the problem
        \begin{align}
            \label{eq:rcprod:relaxed_y}
            \min_{x, t} \quad 
            \left\{
                (d - yr)^{\top}x + f^{\top}t + by
            \ \middle| \
                \eqref{eq:rcprod:conic:primal:conic}
            \right\}.
        \end{align}
        Problem \eqref{eq:rcprod:relaxed_y} is immediately strictly feasible, and bounded since $(d - yr), f \, {>} \, 0$ and $x, t \, {\geq} \, 0$.
        Hence, Theorems \ref{thm:completion:existence} and \ref{thm:completion:optimal} apply, and there exists a dual-optimal completion to recover $\pi, \sigma, \tau$.
        A closed-form completion is then achieved as follows.
        First, constraints \eqref{eq:rcprod:conic:dual:pi} and \eqref{eq:rcprod:conic:dual:tau} yield $\pi = d - ry$ and $\tau = f$.
        Next, note that $\sigma$ only appears in constraint \eqref{eq:rcprod:conic:dual:rsoc} and has negative objective coefficient.
        Further noting that \eqref{eq:rcprod:conic:dual:rsoc} can be written as $\sigma_{j}^{2} \leq 2 \pi_{j} \tau_{j}$, it follows that $\sigma_{j} = - \sqrt{2 \pi_{j} \tau_{j}}$ at the optimum.

    \subsubsection{Numerical Results}
    \label{sec:rcprod:results}

        \begin{table}[!t]
            \centering
            \caption{Comparison of optimality gaps on production planning instances.}
            \label{tab:res:rcprod:gaps}
            \small
            \begin{tabular}{rrrrrrrrr}
                \toprule
                    &
                    & \multicolumn{3}{c}{DC3}
                    & \multicolumn{3}{c}{DLL}\\
                \cmidrule(lr){3-5}
                \cmidrule(lr){6-8}
                \multicolumn{1}{c}{$n$} 
                    & Opt val$^*$
                    & avg. & std & max 
                    & avg. & std & max\\
                \midrule
  10 &   3441.8 &  70.76 &   9.42 &  90.23 &   \textbf{0.23} &   0.57 &  17.05 \\
  20 &   6988.2 &  78.52 &   6.67 &  92.31 &   \textbf{0.41} &   0.69 &   9.04 \\
  50 &  17667.4 &  81.70 &   5.41 &  92.69 &   \textbf{1.03} &   1.69 &  21.68 \\
 100 &  35400.2 &  83.25 &   4.78 &  93.31 &   \textbf{0.37} &   0.57 &   6.69 \\
 200 &  70889.5 &  84.06 &   4.20 &  93.44 &   \textbf{0.29} &   0.46 &   4.81 \\
 500 & 177060.0 &  86.74 &   3.80 &  93.74 &   \textbf{0.46} &   0.73 &   9.92 \\
1000 & 354037.5 &  87.01 &   3.71 &  93.80 &   \textbf{0.36} &   0.48 &   4.44 \\
                \bottomrule
            \end{tabular}\\
            \footnotesize{All gaps are in \%; best values are in bold. $^*$Mean optimal value on test set; obtained with Mosek.}
        \end{table}

        Table \ref{tab:res:rcprod:gaps} reports optimality gap statistics for DC3 and DLL.
        Similar to the linear programming setting, DLL substantially outperforms DC3, with average optimality gaps ranging from 0.23\% to 1.03\%, compared with 70.76\%--87.01\% for DC3.
        In addition, DLL exhibits smaller standard deviation and maximum optimality gaps than DC3.
        These results can be explained by several factors.
        First, the neural network architecture used in DC3 has output size $n+1$, compared to $1$ for DLL; this is because DLL leverages a more efficient dual completion procedure.
        Second, a closer examination of DC3's output reveals that it often fails to satisfy the (conic) inequality constraints \eqref{eq:rcprod:conic:dual:y} and \eqref{eq:rcprod:conic:dual:rsoc}.
        More generally, DC3 was found to have much slower convergence than DLL during training.
        While the performance of DC3 may benefit from more exhaustive hypertuning, doing so comes at a significant computational and environmental cost.
        This further highlights the benefits of DLL, which requires minimal tuning and is efficient to train.

        \begin{table}[!t]
            \centering
            \caption{Computing time statistics for nonlinear instances}
            \label{tab:res:rcprod:timing}
            \small
            \begin{tabular}{rrrr}
                \toprule
                 $n$ & \multicolumn{1}{c}{Mosek$^{\dagger}$} & \multicolumn{1}{c}{DC3$^{\ddagger}$} & \multicolumn{1}{c}{DLL$^{\ddagger}$}  \\
                 \midrule
  10 &     73.5 CPU.s &    2.7 GPU.ms &    0.2 GPU.ms\\
  20 &     75.3 CPU.s &    2.7 GPU.ms &    0.2 GPU.ms\\
  50 &     15.4 CPU.s &    2.7 GPU.ms &    0.2 GPU.ms\\
 100 &     24.9 CPU.s &    2.7 GPU.ms &    0.4 GPU.ms\\
 200 &     49.9 CPU.s &    5.1 GPU.ms &    1.0 GPU.ms\\
 500 &     98.8 CPU.s &   15.9 GPU.ms &    5.1 GPU.ms\\
1000 &    203.0 CPU.s &   41.5 GPU.ms &   19.0 GPU.ms\\
                \bottomrule
            \end{tabular}\\
            {\footnotesize $^{\dagger}$Time to solve all instances in the test set, using one CPU core. $^{\ddagger}$Time to run inference on all instances in the test set, using one V100 GPU.}
        \end{table}

        Finally, Table \ref{tab:res:rcprod:timing} reports computing time statistics for Mosek, a state-of-the-art conic interior-point solver, DC3 and DLL.
        Abnormally high times are observed for Mosek and $n {=} 10, 20$.
        These are most likely caused by congestion on the computing nodes used in the experiments, and are discarded in the analysis.
        Again, DC3 and DLL outperform Mosek by about three orders of magnitude.
        Furthermore, DLL is about 10x faster than DC3 for smaller instances ($n {\leq} 100$), and about 2x faster for the largest instances ($n{=}1000$).
        This is caused by DC3's larger output dimension and correction steps.

\section{Limitations}
\label{sec:limitations}

The main theoretical limitation of the paper is that it considers convex conic optimization problems, and therefore does not consider discrete decisions nor general non-convex constraints.
Since convex relaxations are typically used to solve non-convex problems to global optimality, the proposed approach is nonetheless still useful in non-convex settings.

On the practical side, the optimal dual completion presented in Section \ref{sec:methodology:completion} requires, in general, the use of an implicit layer, which is typically not tractable for large-scale problems.
In the absence of a known closed-form optimal dual completion, it may still be possible to design efficient completion strategies that at least ensure dual feasibility.
One such strategy is to introduce artificial large bounds on all primal variables, and use the completion outlined in Example \ref{ex:bounded_variables}.
Finally, all neural network architectures considered in the experiments are fully-connected neural networks.
Thus, a separate model is trained for each input dimension.
Nevertheless, the DLL methodology is applicable to graph neural network architectures, which would support arbitrary problem size.
The use of GNN models in the DLL context is a promising avenue for future research.

\section{Conclusion}
\label{sec:conclusion}

    The paper has proposed Dual Lagrangian Learning (DLL), a principled methodology for learning dual conic optimization proxies.
    Thereby, a systematic dual conic completion, differentiable conic projection layers, and a self-supervised dual Lagrangian training framework have been proposed.
    The effectiveness of DLL has been demonstrated on numerical experiments that consider linear and nonlinear conic problems, where DLL significantly outperforms DC3 \cite{Donti2021_DC3}, and achieves 1000x speedups over commercial interior-point solvers.

    One of the main advantages of DLL is its simplicity.
    The proposed dual completion can be stated only in terms of \emph{primal} constraints, thus relieving users from the need to explicitly write the dual problem.
    DLL introduces very few hyper-parameters, and requires minimal tuning to achieve good performance.
    This results in simpler models and improved performance, thus delivering computational and environmental benefits.

    DLL opens the door to multiple avenues for future research, at the intersection of ML and optimization.
    The availability of high-quality dual-feasible solutions naturally calls for the integration of DLL in existing optimization algorithms, either as a warm-start, or to obtain good dual bounds fast.
    Multiple optimization algorithms have been proposed to optimize Lagrangian functions, which may yield more efficient training algorithms in DLL.
    Finally, given the importance of conic optimization in numerous real-life applications, DLL can provide a useful complement to existing primal proxies.

\section*{Acknowledgements}

    This research was partially supported by NSF awards 2007164 and 2112533, and ARPA-E PERFORM award DE-AR0001280.

\newpage
\bibliographystyle{unsrtnat}
\bibliography{refs}

\clearpage
\appendix
\section{Standard cones}
\label{app:cones}

    This section presents standard cones and their duals, as well as corresponding Euclidean and radial projections.
    The reader is referred to \cite{Coey2022_Hypatia} for a more exhaustive list of non-standard cones, and to \cite[Sec. 6.3]{Parikh2014_ProximalAlgorithms} for an overview of Euclidean projections onto standard cones.

    \subsection{Non-negative orthant}
    \label{app:cones:R+}
    
        The non-negative orthant is defined as $\R^{n}_{+} = \{x \in \R^{n}: x \geq 0\}$.
        It is a self-dual cone, and forms the basis of linear programming \cite{ben2001lectures}.

        \paragraph{Euclidean projection}
        The Euclidean projection on $\R^{n}_{+}$ is
        \begin{align}
            \label{eq:cones:R+:proj:euclidean}
            \Pi_{\R^{n}_{+}}(\ybar) = \max(0, \ybar) = \relu(\ybar),
        \end{align}
        where the $\max$ and $\relu$ operations are performed element-wise.

        \paragraph{Radial projection}
        The radial projection with ray $\e$, applied coordinate-wise, is equivalent to the Euclidean projection.

    \subsection{Conic quadratic cones}
    \label{app:cones:SOC}

        Conic quadratic cones include the second-order cone (SOC)
        \begin{align}
            \label{eq:cones:SOC}
            \Q^{n} = \{
                x \in \R^{n} 
                : 
                x_{1} \geq \sqrt{x^{2}_{2} {+} \cdots {+} x_{n}^{2}}
            \}
        \end{align}
        and the rotated second-order cone (RSOC)
        \begin{align}
            \label{eq:cones:RSOC}
            \Qr^{n} = \{ x \, {\in} \R^{n} \, {:} \, 2 x_{1} x_{2} \, {\geq} \,
            x^{2}_{3} {+} \cdots {+} x_{n}^{2}, x_{1}, x_{2} \geq 0
            \}.
        \end{align}
        Both cones are self-dual, i.e., $\Q^{*} = \Q$ and $\Qr^{*} = \Qr$.
        The RSOC is the main building block of conic formulations of convex quadratically-constrained optimization problems.

        \paragraph{Euclidean projection}
        The Euclidean projection on $\Q^{n}$ is given by
        \begin{align}
            \label{eq:cones:SOC:proj:euclidean}
            \Pi_{\Q^{n}}(\xbar) = 
            \left\{
            \begin{array}{ll}
                \xbar & \text{ if } \xbar \in \Q^{n}\\
                0 & \text{ if } \xbar \in -\Q^{n}\\
                \frac{\xbar_{1} + \delta}{2 \delta} (\delta, \xbar_{2}, ..., \xbar_{n}) & \text{ otherwise}
            \end{array}
            \right.
        \end{align}
        where $\delta \, {=} \, \|(\xbar_{2}, ..., \xbar_{n})\|_{2}$.

        \paragraph{Radial projection}
        Given interior ray $\rho = (1, 0, ..., 0) \kge_{\Q^{n}} 0$, the radial projection is
        \begin{align}
            \label{eq:cones:SOC:proj:radial}
            \Pi_{\Q^{n}}^{\rho}(\xbar) &= (\xhat_{1}, \xbar_{2}, ..., \xbar_{n}),
            &
            \xhat_{1} &= \max(\xbar_{1}, \|(\xbar_{2}, ..., \xbar_{n})\|_{2}).
        \end{align}
        Note that, in the worst case, computing $\Pi_{\Q}(\xbar)$ requires $\mathcal{O}(2n)$ operations, and modifies all coordinates of $\bar{x}$.
        In contrast, computing $\Pi_{\Q}^{\rho}(\xbar)$ requires only $\mathcal{O}(n)$ operations, and only modifies the first coordinate of $\xbar$.

        Closed-form formulae for Euclidean and radial projections onto $\Qr^{n}$ are derived from \eqref{eq:cones:SOC:proj:euclidean} and \eqref{eq:cones:SOC:proj:radial}.

    \subsection{Positive Semi-Definite cone}
    \label{app:cones:SDP}
    
        The cone of positive semi-definite (PSD) matrices of order $n$ is defined as
        \begin{align}
            \Ksdp^{n} &= \{ X \in \R^{n\times n}: X = X^{\top}, \eigmin(X) \geq 0 \}.
        \end{align}
        Note that all matrices in $\Ksdp^{n}$ are symmetric, hence all their eigenvalues are real.
        The PSD cone is self-dual, and generalizes the non-negative orthant and SOC cones \cite{boyd2004convex}.

        \paragraph{Euclidean projection}
        The Euclidean projection onto $\Ksdp^{n}$ is given by
        \begin{align}
            \label{eq:cones:SDP:proj:euclidean}
            \Pi_{\Ksdp^{n}}(\bar{X}) &= \sum_{i} \max(0, \lambda_{i}) v_{i} v_{i}^{\top},
        \end{align}
        where $\bar{X} \in \R^{n \times n}$ is symmetric with eigenvalue decomposition $\bar{X} = \sum_{i} \lambda_{i} v_{i} v_{i}^{\top}$.
        Note that the Euclidean projection onto the PSD code thus requires a full eigenvalue decomposition, which has complexity $\mathcal{O}(n^{3})$.

        \paragraph{Radial projection}
        The radial projection considered in the paper uses $\rho = I_{n} \in \intr (\Ksdp^{n})$.
        This yields the closed-form projection
        \begin{align}
            \label{eq:cones:SDP:proj:radial}
            \Pi_{\Ksdp^{n}}^{\rho}(\bar{X}) &= \bar{X} + \min(0, |\eigmin(\bar{X})|) I_{n}.
        \end{align}
        Note that the radial projection only requires computing the smallest eigenvalue of $\bar{X}$, which is typically much faster than a full eigenvalue decomposition, and only modifies the diagonal of $\bar{X}$.

    \subsection{Exponential Cone}
    \label{app:cones:exp}
    
        The 3-dimensional exponential cone is a non-symmetric cone defined as
        \begin{align}
            \Kexp 
            &= 
            \cl \left\{
                x \in \R^{3}
            :
                x_{1} \geq x_{2} e^{x_{3} / x_{2}},
                x_{2} > 0
            \right\},
        \end{align}
        whose dual cone is
        \begin{align}
            \Kexp^{*}
            &= \cl
            \left\{
                y \in \R^{3}
            :
                \frac{-y_{1}}{y_{3}} \geq e^{\frac{y_{2}}{y_{3}} - 1},
                y_{1} > 0,
                y_{3} < 0
            \right\}.
        \end{align}
        The exponential cone is useful to model exponential and logarithmic terms, which occur in, e.g., relative entropy, logistic regression, or logarithmic utility functions.

        \paragraph{Euclidean projection}
        To the best of the authors' knowledge, there is no closed-form, analytical formula for evaluating $\Pi_{\Kexp}$ nor $\Pi_{\Kexp^{*}}$, which instead require a numerical method, see, e.g., \cite{Parikh2014_ProximalAlgorithms} and \cite{Friberg2023_ProjectionExponentialCone} for completeness.

        \paragraph{Radial projection}
        To avoid any root-finding operation, the paper leverages the fact that $x_{1}, x_{2} \, {>} \, 0, \forall (x_{1}, x_{2}, x_{3}) \in \Kexp$.
        Note that one can enforce $\bar{x}_{1}, \bar{x}_{2} \, {>} \, 0$ via, e.g., softplus activation.
        A radial projection is then obtained using $\rho = (0, 0, 1)$, which yields
        \begin{align}
            \label{eq:cones:exp:proj:radial:primal}
            \Pi_{\Kexp}^{\rho}(\bar{x}_{1}, \bar{x}_{2}, \bar{x}_{3})
            &= \left(
                \bar{x}_{1},
                \bar{x}_{2},
                \min \big( \bar{x}_{3}, \bar{x}_{2} \log \frac{\bar{x}_{1}}{\bar{x}_{2}} \big)
            \right).
        \end{align}
        This approach does not require any root-finding, and is therefore more amenable to automatic differentiation.
        The validity of Eq. \eqref{eq:cones:exp:proj:radial:primal} is immediate from the representation
        \begin{align}
            \Kexp = \cl \{x \in \R^{3}: \frac{x_{3}}{x_{2}} \leq  \log \frac{x_{1}}{x_{2}}, x_{1}, x_{2} > 0\}.
        \end{align}
        Similarly, assuming  $\ybar_{1} > 0$ and $\ybar_{3} < 0$, the radial projection onto $\Kexp^{*}$ reads
        \begin{align}
            \label{eq:cones:exp:proj:radial:dual}
            \Pi_{\Kexp^{*}}^{\rho}(\bar{x}_{1}, \bar{x}_{2}, \bar{x}_{3})
            &= \left(
                \bar{y}_{1},
                \max \big(\ybar_{2}, \ybar_{3} {+} \ybar_{3} \ln \frac{\ybar_{1}}{-\ybar_{3}}\big),
                \bar{y}_{3}
            \right).
        \end{align}

    \subsection{Power Cone}
    \label{app:cones:pow}
    
        Given $0 \, {<} \, \alpha \, {<} \, 1$, the 3-dimensional power cone is defined as
        \begin{align}
            \Kpow_{\alpha}
            &= 
            \left\{
                x \in \R^{3}
            :
                x_{1}^{\alpha} x_{2}^{1-\alpha}
                \geq
                |x_{3}|,
                x_{1}, x_{2} \geq 0
            \right\}.
        \end{align}
        Power cones are non-symmetric cones, which allow to express power other than $2$, e.g., $p$-norms with $p \geq 1$.
        Note that $\Kpow_{\nicefrac{1}{2}}$ is a scaled version of the rotated second-order cone $\Qr^{3}$.
        The 3-dimensional power cone $\Kpow_{\alpha}$ is sufficient to express more general, high-dimensional power cones.
        The dual power cone is
        \begin{align}
            \Kpow_{\alpha}^{*}
            &= 
            \left\{
                y \in \R^{3}
            :
                \left(
                    \frac{y_{1}}{\alpha}, \frac{y_{2}}{1-\alpha}, y_{3}
                \right)
                \in \Kpow_{\alpha}
            \right\}.
        \end{align}

        \paragraph{Euclidean projection}
        The Euclidean projection onto the power cone $\Kpow_{\alpha}$ is described in \cite{Hien2015_EuclideanProjectionPowerCone}.
        Similar to the exponential cone, it requires a root-finding operation.

        \paragraph{Radial projection}
        The proposed radial projection is similar to the one proposed for $\Kexp$.
        Assuming $\xbar_{2}, \xbar_{3} \, {>} \, 0$, and using $\rho = (1, 0, 0)$, the radial projection reads
        \begin{align}
            \label{eq:cones:pow:proj:radial}
            \Pi_{\Kpow_{\alpha}}^{\rho}(\bar{x}_{1}, \bar{x}_{2}, \bar{x}_{3})
            &= \left(
                \max \big(
                    \bar{x}_{1}, 
                    \bar{x}_{2}^{\frac{\alpha-1}{\alpha}} |\bar{x}_{3}|^{\frac{1}{\alpha}}
                \big),
                \bar{x}_{2},
                \bar{x}_{3}
            \right).
        \end{align}
        A similar approach is done to recover $y \in \Kpow_{\alpha}^{*}$ after scaling the first two coordinates of $y$.
        This technique can be extended to the more general power cones.

\section{Experiment Details}
\label{sec:app:experiment_details}

    \subsection{Common experiment settings}
    \label{sec:app:experiment_details:common}
    
        All experiments are conducted on the Phoenix cluster \cite{PACE}, on Linux machines equipped with Intel Xeon Gold 6226 CPU @ 2.70GHz nodes and using a single Tesla V100 GPU.
        Each job was allocated 12 CPU threads and 64GB of RAM.
        All ML models are formulated and trained using Flux \cite{FluxML}; unless specified otherwise, all (sub)gradients are computed using the auto-differentiation backend Zygote \cite{ZygoteAD}.
        All linear problems are solved with Gurobi v10 \cite{Gurobi}.
        All nonlinear conic problems are solved with Mosek \cite{Mosek}.

        All neural network architectures considered here are fully-connected neural networks (FCNNs).
        Thus, a separate model is trained for each input dimension.
        Note that the proposed DLL methodology is applicable to graph neural network architectures, which would support arbitrary problem size.
        The use GNN models in the DLL context is a promising avenue for future research; a systematic comparison of the performance of GNN and FCNN architectures is, however, beyond the scope of this work.

        All ML models are trained in a self-supervised fashion following the training scheme outlined in Section \ref{sec:methodology:training}, and training is performed using the Adam optimizer \cite{kingma2014_Adam}.
        The training scheme uses a patience mechanism where the learning rate $\eta$ is decreased by a factor $2$ if the validation loss does not improve for more than $N_{p}$ epochs.
        The initial learning rate is $\eta \, {=} \, 10^{-4}$.
        Training is stopped if either the learning rate reaches $\eta_{\text{min}} \, {=} \, 10^{-7}$, or a maximum $N_{\text{max}}$ epochs is reached.
        Every ML model considered in the experiments was trained in under an hour.

        A limited, manual, hypertuning was performed by the authors during preliminary experiments.
        It was found that DLL models require very little hypertuning, if any, to achieve satisfactory performance.
        In contrast, DC3 was found to require very careful hypertuning, even just to ensure its numerical stability.
        It is also important to note that DC3 introduces multiple additional hyperparameters, such as the number of correction steps, learning rate for the correction steps, penalty coefficient for the soft penalty loss, etc.
        These additional hyperparameters complicate the hypertuning task, and result in additional computational needs.
        Given the corresponding economical and environmental costs, only limited hypertuning of DC3 was performed.

        Finally, it was observed that DC3 often fail to output dual-feasible solutions, which therefore do not valid dual bounds.
        Therefore, to ensure a fair comparison, the dual solution produced by DC3 is fed to the dual optimal completion of DLL, thus ensuring dual feasibility and a valid dual bound.
        This is only performed at test time, with a negligible overhead since the dual completion uses a closed-form formula.
        All optimality gaps for DC3 are reported for this valid dual bound.

    \subsection{Linear programming problems}
    \label{sec:app:experiment_details:MDK}

        \subsubsection{Problem formulation}
            The first set of experiments considers the continuous relaxation of multi-dimensional knapsack problems \cite{Freville1994_MDKnapsackPreprocessing,Freville2004_MultiDimensionalKnapsackOverview}, which are of the form
            \begin{align}
                \min_{x} \quad
                \left\{
                    -p^{\top}x
                \ \middle| \ 
                    Wx \leq b,
                    x \in [0, 1]^{n}
                \right\},
            \end{align}
            where $n$ denotes the number of items, $m$ denotes the number of resources, $p \in \R_{+}^{n}$ is the value of each item, $b_{i}$ is the amount of resource $i$, and $W_{ij}$ denotes the amount of resource $i$ used by item $j$.
            The dual problem reads
            \begin{align}
                \max_{y, \zl, \zu} \quad
                \left\{
                    b^{\top} y - \e^{\top} \zu
                \ \middle| \
                    W^{\top}y + \zl - \zu = -p,
                    y \leq 0, \zl \geq 0, \zu \geq 0,
                \right\}
            \end{align}

        \subsubsection{Data generation}
        
            For each number of items $n \, {\in} \, \{100, 200, 500\}$ and number of resources $m \, {\in} \, \{5, 10, 30\}$, a total of 
            $16384$ instances are generated using the same procedure as the MIPLearn library \cite{SantosXavier2023_MIPLearn}.
            Each instance is solved with Gurobi, and the optimal dual solution is recorded for evaluation purposes.
            This dataset is split in training, validation and testing sets, which contain 8192, 4096 and 4096 instances, respectively.

        \subsubsection{DLL implementation}
    
            The DLL architecture considered here is a fully-connected neural network (FCNN); a separate model is trained for each combination $(m, n)$.
            The FCNN model takes as input the flatted problem data $(b, p, W) \in \mathbb{R}^{1+n+n \times m}$, and outputs $y \in \R^{m}$.
            The FCNN has two hidden layers of size $2(m+n)$ and sigmoid activation; the output layer uses a negated softplus activation to ensure $y \leq 0$.
            The dual completion procedure follows Example \eqref{ex:bounded_variables}.

            \paragraph{Hyperparameters}
                The patience parameter is $N_{p} \, {=} \, 32$, and the maximum number of training epochs is $N_{\text{max}} \, {=} \, 1024$.

        \subsubsection{DC3 implementation}

            The DC3 architecture consists of an initial FCNN which takes as input $(b, p, W)$, and outputs $y, \zl$.
            Then, $\zu$ is recovered by equality completion as $\zu = p + W^{\top}y + \zl$.
            The correction then applies gradient steps $(y, \zl) \gets (y, \zl) - \gamma \nabla \phi(y, \zl)$ where
            \begin{align*}
                \phi(y, \zl) = \|\max(0, y) \|^{2} + \|\min(0, \zl)\|^{2} + \| \min(0, \zu) \|^{2}
            \end{align*}
            The corresponding gradients $\nabla \phi(y, \zl)$ were computed analytically.
            After applying corrections, the dual equality completion is applied one more time to recover $\zu$, and the final soft loss is
            \begin{align}
                b^{\top}y - \e^{\top} \zu + \rho \times \phi(y, \zl)
            \end{align}
            which considers both the dual objective value, and the violation of inequality constraints.
            Note that the dual objective $b^{\top}y - \e^{\top} \zu$ is not a valid dual bound in general, because $y, \zl, \zu$ may not be dual-feasible.

            \paragraph{Hyperparameters}
                The maximum number of correction steps is $10$, the learning rate for correction is $\gamma \, {=} \, 10^{-4}$.
                The soft penalty weight is set to $\rho = 10$; this parameter was found to have a high impact on the numerical stability of training.
                The patience parameter is $N_{p} \, {=} \, 32$, and the maximum number of training epochs is $N_{\text{max}} \, {=} \, 1024$.

    \subsection{Nonlinear Production and Inventory Planning Problems}
    \label{sec:app:experiment_details:rcprod}

        \subsubsection{Problem formulation}
        
            The original presentation of the resource-constrained production and inventory planning problem  \cite{Ziegler1982_ProductionPlanning} uses the nonlinear convex formulation
            \begin{subequations}
            \label{eq:rcprod:NLP}
            \begin{align}
                \min_{x} \quad
                & \sum_{j} d_{j} x_{j} + f_{j} \frac{1}{x_{j}}\\
                s.t. \quad
                & r^{\top}x \leq b,\\
                & x \geq 0,
            \end{align}
            \end{subequations}
            where $n$ is the number of items to be produced, $x \in \R^{n}$, $b \in \R$ denotes the available resource amount, and $r_{j} > 0$ denotes the resource consumption rate of item $j$.
            The objective function captures production and inventory costs.
            Namely, $d_{j} \, {=} \, \frac{1}{2} c_{j}^{p} c_{j}^{r}$ and $f_{j} \, {=} \, c^{o}_{j} D_{j}$, where $c^{p}_{j}, c^{r}_{j}, c^{o}_{j} \, {>} \, 0$ and $D_{j} \, {>} \, 0$ denote per-unit holding cost, rate of holding cost, ordering cost, and total demand for item $j$, respectively.
            
            The problem is reformulated in conic form \cite{MosekModelingCookbook} as
            \begin{subequations}
            \begin{align}
                \min_{x, t} \quad
                & d^{\top}x + f^{\top}t\\
                s.t. \quad
                & r^{\top}x \leq b,\\
                & (x_{j}, t_{j}, \sqrt{2}) \in \Qr^{3}, \quad \forall j=1, ..., n.
            \end{align}
            \end{subequations}
            whose dual problem is
            \begin{subequations}
                \begin{align}
                    \max_{y, \pi, \tau, \sigma} \quad
                    & b y - \sqrt{2} \e^{\top} \sigma_{j}\\
                    s.t. \quad 
                    & r y + \pi = d,\\
                    & \tau = f,\\
                    & y \leq 0,\\
                    & (\pi_{j}, \tau_{j}, \sigma_{j}) \in \Qr^{3},  \quad \forall j=1...n.
                \end{align}
            \end{subequations}
            Note that the dual problem contains $1+3n$ variables, $2n$ equality constraints, $1$ linear inequality constraints, and $n$ conic inequality constraints.
            Therefore, DC3 must predict $n+1$ dual variables, then recover $2n$ variables by equality completion, and correct for $n+1$ inequality constraints.
            In contrast, by exploiting dual optimality conditions, DLL eliminates $3n$ dual variables, thus reducing the output dimension of the initial prediction from $n+1$ to $1$, and eliminates the need for correction.

        \subsubsection{Data generation}

            For each $n \in \{10, 20, 50, 100, 200, 500, 1000\}$, 16384 instances are generated using the procedure of \cite{Ziegler1982_ProductionPlanning}.
            First, $D_{j}$ is sampled from a uniform distribution $U[1, 100]$, $c^{p}_{j}$ is sampled from $U[1, 10]$, and $c^{r}_{j}$ is sampled from $U[0.05, 0.2]$.
            Then, $c^{o}_{j} = \alpha_{j} c^{p}_{j}$ and $r_{j} = \beta_{j} c^{p}_{j}$, where $\alpha, \beta$ are sampled from $U[0.1, 1.5]$ and $U[0.1, 2]$, respectively.
            Finally, the right-hand side is $b \, {=} \, \eta \sum_{j} r_{j}$ where $\eta$ is sampled from $U[0.25, 0.75]$.
    
            Each instance is solved with Mosek, and its solution is recorded for evaluation purposes.
            The dataset is split into training, validation and testing sets comprising 8192, 4096 and 4096 instances, respectively.

        \subsubsection{DLL implementation}

            The DLL architecture consists of an initial FCNN that takes as input $(d, f, r, b) \, {\in} \, \R^{1+3n}$, and output $y \in \mathbb{R}$. 
            The FCNNs have two hidden layers of size $\max(128, 4n)$ and sigmoid activation.
            For the output layer, a negated softplus activation ensures $y \leq 0$.
            The dual completion outlined in Section \ref{sec:exp:rcprod} then recovers $(\pi, \sigma, \tau)$.

            \paragraph{Hyperparameters}
                The patience parameter is $N_{p} \, {=} \, 128$, and the maximum number of training epochs is $N_{\text{max}} \, {=} \, 4096$.
                The patience mechanism is deactivated for the first 1024 epochs; this latter setting has little impact of the performance of DLL, and was introduced to avoid premature termination for DC3.

        \subsubsection{DC3 implementation}

            The DLL architecture consists of an initial FCNN that takes as input $(d, f, r, b) \, {\in} \, \R^{1+3n}$, and outputs $y, \sigma \in \mathbb{R}$. 
            The FCNNs have two hidden layers of size $\max(128, 4n)$ and sigmoid activation, and the output layer has linear activation.

            The equality completion step recovers $\pi = d - ry$ and $\tau = f$.
            The correction step then apply gradient steps to minimize the violations $\phi(y, \sigma) = \phi_{y}(y, \sigma) + \phi_{\pi}(y, \sigma) + \phi_{\sigma}(y, \sigma)$, where
            \begin{align}
                \phi_{y}(y, \sigma)      &= \max(0, y)^{2},\\
                \phi_{\pi}(y, \sigma)    &= \min(0, \pi)^{2},\\
                \phi_{\sigma}(y, \sigma) &= \sum_{j} \max(0, \sigma_{j}^{2} - 2 \pi_{j} \tau_{j})^{2}.
            \end{align}
            Note that, to express $\phi_{\sigma}(y, \sigma)$, conic constraints \eqref{eq:rcprod:conic:dual:rsoc} are converted to their nonlinear programming equality, because DC3 does not handle conic constraints.
            Gradients for $\phi$ are computed analytically, and implemented directly in the inequality correction procedure.
            The final soft loss is then            
            \begin{align}
                b y - \sqrt{2} \e^{\top} \sigma_{j} + \rho \phi(y, \sigma).
            \end{align}
            
            \paragraph{Hyperparameters}
                The maximum number of correction steps is $10$, the learning rate for correction is $\gamma \, {=} \, 10^{-5}$, and the soft loss penalty parameter is $\rho=10$.
                The patience parameter is $N_{p} \, {=} \, 128$, and the maximum number of training epochs is $N_{\text{max}} \, {=} \, 4096$.
                The patience mechanism is deactivated for the first 1024 epochs; this latter setting has little impact of the performance of DLL, and was introduced to avoid premature termination for DC3.

\end{document}